\documentclass{amsart}

\usepackage{amsmath}
\usepackage{amssymb}
\usepackage{amsfonts}
\usepackage{amsthm}
\usepackage{float}

\newtheorem{teo}{Theorem}[section]
\newtheorem{coro}[teo]{Corollary}
\newtheorem{lem}[teo]{Lemma}
\theoremstyle{definition}
\newtheorem{defi}{Definition}[section]
\newtheorem{ejems}[defi]{Examples}
\newtheorem{ejem}[defi]{Example}

\newcommand{\R}{\mathbb{R}}
\newcommand{\N}{\mathbb{N}} 

\newcommand{\tlim}{\displaystyle\lim}

\begin{document}

\title[Fixed points of $D_T(a,b)$]{\large Fixed point theorems for a class of mappings depending of another function
 and defined on cone metric spaces}
\author[J. R. Morales and E.M. Rojas]{Jos\'e R. Morales and Edixon Rojas\\\\
\Small University Of Los Andes, Faculty Of Science\\
Department Of Mathematics, 5101\\
 M\'Erida, Venezuela \\
e-mail: moralesj@ula.ve, edixonr@ula.ve\\
}

\date{\today}
 \maketitle

\begin{abstract}
In this paper we study the existence and uniqueness of fixed points of a class of mappings defined on complete, (sequentially compact) cone metric spaces, without continuity conditions and depending on another function.
\end{abstract}

%

\section{Introduction}
Since the Banach's contraction principle come out in 1922 (see, e.g., \cite{KK}, \cite{8}), several type of contraction mappings, acting on metric spaces, have been appearing, at point that in 1977 B.E. Rhoades \cite{Roa77} made a comparison of 250 different type (or classes) of contraction mappings. Chatterjea, Eldestein, Hardy \& Rogers,  Kannan and Shimi -contraction type of maps, are some of the most extensively studied \cite{1,6,7,9,91,14,15}. Generalizations of  such classes are already known, we quote here only \cite{Be93}, \cite{Ci74}, \cite{DaDe86}, \cite{RV} and \cite{RaCh94} for some of such generalizations.  Also, there exist a lot of  works devoted to study the relations among these (and others) generalized contraction-type maps, see for instance \cite{CoCar97,Roa83,Roa88}.
Recently, a new generalization of contraction mappings acting on complete metric spaces was introduced by A. Beiranvand, S. Moradi, M. Omid and H. Pazandeh \cite{BMOP} called: $T-$contraction and $T-$contrative which are contraction mappings depending on another function $T$, extending in this fashion the Banach's contraction principle and the Edelstein's fixed point Theorem. Furthermore, in \cite{Mo09} S. Moradi  define the corresponding $T$-Kannan contraction.

The class of contraction mappings that will be play a principal role in this paper is the  introduced in 1986 by
L. Nova \cite{10}, (or see \cite{2}, \cite{3} for instance) which is defined as follows:

\emph{Let $(M,d)$ be a complete metric space, $K\subset M$ closed and $T:K\longrightarrow K$ an arbitrary mapping.
We shall say $T$ belong or is of class $D(a,b)$ when $T$ satisfies the following condition, for $a,b\geq0$ and any $x,y\in K$.
\begin{equation}\label{Dab}
d(Tx,Ty)\leq ad(x,y)+b\left[d(x,Tx)+d(y,Ty)\right].
\end{equation}}
Fixed point's existence and uniqueness theorems, as well  examples of mappings on $D(a,b)$, can be found in  \cite{2,3,10,10Y}.
 Notice that a mapping $T$ satisfying \eqref{Dab} is a Banach contraction map when $b=0$ and $a<1$. In addition, for $a=0$ we obtain that $T$ is a Kannan contraction and  for  $a\neq0$, a Reich contraction for the special case of $a_2=a_3$, while any mapping, continuous or not, is in class $D(1,1)$.

On the other hand, generalized Banach spaces, specially ordered Banach spaces, bring a large variety of spaces where contraction type of mappings can be defined, providing conditions for the existence and uniqueness of fixed points (see,  \cite{Be91,Bra00,Mo09,Z}).  Using this approach, Huan Long - Guang and Zhang Xian in \cite{GX}, introduced the concept of cone metric space, where the set of real numbers is replaced by an ordered Banach space.  They introduced the basic definitions and discuss some properties of convergence of sequences in cone metric spaces. Ordered Banach spaces and normal cones are especially useful in applications on optimization theory coming from non-convex analysis, see for instance \cite{HM,MSR}.

In the same paper \cite{GX}, several fixed point theorems for contractive single-valued maps (Banach, Eldestein and Kannan type-mappings) in such spaces was proved. Subsequently, some other mathematicians,  (see, e.g., \cite{AJ,AR,DR,JR,RV,RH}), have generalized these results of Guang and Zhang. The corresponding version of $D(a,b)$ defined on cone metric spaces was given in 2008 by M.S. Khan and M. Samanipour  \cite{KhSa08} as a particular case of the class $E(a,b,c)$ introduced there (Definition 2.4).

The main goal of this paper is, using the approach of A. Beiranvand \emph{et al} \cite{BMOP,Mo09} and M.S. Khan and M. Samanipour \cite{KhSa08}, to introduce a generalization of the class $D(a,b)$ and proof the corresponding theorems of existence and uniqueness of fixed points.

\section{Cone metric spaces}
Consistent with Guang and Zhang \cite{GX}, we recall the definitions of cone metric space, the notion of convergence and other results that will be needed in the sequel.

Let $E$ be a real Banach space and $P$ a subset of $E.\,\, P$ is called a \emph{cone} if and only if:
\begin{itemize}
\item[(1)] $P$ is non-empty, closed and $P\neq \{0\};$
\item[(2)] $a,b\in \R,\,\, a,b\geq 0$ and $x,y\in P\Rightarrow ax+by\in P;$
\item[(3)] $x\in P$ and $-x\in P\Rightarrow x=0$.
\end{itemize}
For a given cone $P\subseteq E,$ we can define a partial ordering $\leq$ on $E$ with respect to $P$ by
\begin{equation*}
x\leq y,\quad \mbox{if and only if}\quad y-x\in P.
\end{equation*}
We shall write $x<y$ to indicate that $x\leq y$ but $x\neq y,$ while $x\ll y$ will stands for $y-x\in \mbox{Int}\, P,$ where $\mbox{Int}\, P$ denotes the interior of $P.$ The cone $P\subset E$ is called \emph{normal} if there is a number $K>0$ such that for all $x,y\in E,$
$$0\leq x\leq y,\quad \mbox{implies}\quad \|x\|\leq K\|y\|.$$
The least positive number satisfying inequality above is called the \emph{normal constant} of $P$.

The cone $P$ is called \emph{regular} if every increasing sequence which is bounded from above is convergent. That is, if $(x_n)$ is a sequence such that
$$x_1\leq x_2\leq \ldots \leq x_n\leq \ldots \leq y$$
for some $y\in E,$ then there is $x\in E$ such that $\|x_n-x\|\longrightarrow 0,\,\, (n\rightarrow \infty).$

In \cite{RH} was proved that every regular cone is normal, also other properties of the cones can be found there.

In the following we always suppose $E$ is a real Banach space, $P$ is a cone with $\mbox{Int}\, P\neq \emptyset$ and $\leq$ is a partial ordering with respect to $P$.
\begin{defi}[\cite{GX}]
Let $M$ be a nonempty set. Suppose the mapping $d: M\times M\longrightarrow E$ satisfies:
\begin{itemize}
\item[(1)] $0<d(x,y)$ for all $x,y\in M$ and $d(x,y)=0$ if and only if $x=y;$
\item[(2)] $d(x,y)=d(y,x)$ for all $x,y\in M;$
\item[(3)] $d(x,y)\leq d(x,z)+d(y,z)$ for all $x,y,z\in M.$
\end{itemize}
Then $d$ is called a \emph{cone metric} on $M$ and $(M,d)$ is called a \emph{cone metric space}.
\end{defi}
It is obvious that cone metric spaces generalize metric spaces.

\begin{ejems}
\begin{enumerate}
\item \textup{(\cite[Example 1]{GX})} Let $E=\R^2$, $P=\{(x,y)\in E\,:\, x,y\geq 0\}\subset \R^2,\,\, M=\R$ and $d: M\times
      M\longrightarrow E$ such that
\begin{equation*}
 d(x,y)=\bigg(|x-y|,\,\, \alpha|x-y|\bigg)
\end{equation*}
       where $\alpha\geq 0$ is a constant. Then $(M,d)$ is a cone metric space.
\item \textup{(\cite[Example 2.3]{RH})} Let $E=(C^2[0,1],\R)$, with the norm
\begin{equation*}
\|f\|:=\|f\|_\infty+\|f'\|_\infty,
\end{equation*}
and the cone $P=\{\varphi\in E\,:\, \varphi\geq 0\}$. This is an example of a non-normal cone.
\end{enumerate}
\end{ejems}

\begin{defi}[\cite{GX}]
Let $(M,d)$ be a cone metric space and $(x_n)$ a sequence in $M$. Then:
\begin{itemize}
\item[(1)] $(x_n)$ converges to $x\in M$ if for every $c\in E$, with $0\ll c$, there is $n_0\in \N$ such
           that for all $n\geq n_0,$

           $$d(x_n,x)\ll c.$$ We denote this by $\tlim_{n\rightarrow \infty} x_n=x$ or $x_n\longrightarrow x,\,\, (n\rightarrow \infty).$
\item[(2)] If for any $c\in E,$ there is a number $n_0\in \N$ such that for all $m,n\geq n_0$

            $$d(x_n,x_m)\ll c,$$ then $(x_n)$ is called a \emph{Cauchy sequence} in $M;$
\item[(3)] $(M,d)$ is a \emph{complete cone metric space} if every Cauchy sequence is convergent in $M.$
\end{itemize}
\end{defi}

\begin{defi}
Let $(M,d)$ be a cone metric space. If for any sequence $(x_n)$ in $M,$ there is a subsequence $(x_{n_i})$ of $(x_n)$ such that $(x_{n_i})$ is convergent in $M$, then $(M,d)$ is called a \emph{sequentially compact cone metric space}.
\end{defi}

Next Definition and subsequent Lemma are given in \cite{BMOP} in the scope of metric spaces, here we will rewrite it in terms of cone metric spaces.
\begin{defi}
Let $(M,d)$ be a cone metric space. Then
\begin{itemize}
\item[(1)] $T$ is said to be \emph{continuous} if $\tlim_{n\rightarrow \infty} x_n=x$, implies that
           $\tlim_{n\rightarrow \infty} Tx_n=Tx$ for every $(x_n)$ in $M;$

\item[(2)] $T$ is said to be \emph{sequentially convergent} if we have, for every sequence $(y_n),$ if
            $T(y_n)$ is convergent, then $(y_n)$ is also convergent;

\item[(3)] $T$ is said to be \emph{subsequentially convergent} if we have, for every sequence $(y_n),$ if
             $T(y_n)$ is convergent, then $(y_n)$ has a convergent subsequence.
\end{itemize}
\end{defi}

\begin{lem}
Let $(M,d)$ be a sequentially compact cone metric space. Then every function $T: M\longrightarrow M$ is subsequentially convergent and every continuous function $T: M\longrightarrow M$ is sequentially convergent.
\end{lem}

\section{The class $D_T(a,b)$ }

In this section, we will introduce a class of mappings that generalize the class $D(a,b)$ defined in \cite{10} and \cite{KhSa08}, for Banach spaces and complete metric cone spaces respectively.  As in  \cite{RH} our fixed point theorems for elements in this class will be proved omitting the assumption of  normal cone.

\begin{defi}
Let $(M,d)$ be a complete cone metric space and $T,S: M\longrightarrow M$ two functions such that $T$ is continuous and injective. A mapping $S$ is of class $D_T(a,b)$, for $a,b\geq0$ constants, if satisfies the following inequality
\begin{equation}\label{eq3.1}
d(TSx,TSy)\leq a d(Tx,Ty)+b[d(Tx,TSx)+d(Ty,TSy)]
\end{equation} for all $x,y\in M$.
\end{defi}

\begin{ejem}
We are going to  consider a function $S\notin D(a,b)$ and a function $T$, injective and continuous, such that $S\notin D_T(a,b)$.

Let $E=(C[0,1],\R)$, the cone $P=\{\varphi\in E\;:\;\varphi\geq0\}$,  $M:=[0,1]$ and $d(x,y)=|x-y|e^t$ where $e^t\in E$. Let consider $Sx=\sqrt{x}$, the continuous and injective function $Tx=\alpha x$, $\alpha\in\R\setminus\{0\}$. First we will prove that for $0\leq a,b\leq1$,  $S\notin D(a,b)$. In fact, if we suppose the opposite, then
\begin{equation*}
d(Sx,Sy)\leq ad(x,y)+b[d(x,Tx)+d(y,Ty)].
\end{equation*}
Rewriting this inequality we obtain
\begin{equation*}
|Sx-Sy|e^t\leq a|x-y|e^t+b[|x-Sx|+|y-Sy|]e^t,
\end{equation*}
from last inequality we conclude that $S\in D(a,b)$ in the sense of Banach spaces, which means that $S\in D(a,b)$ for the Banach space $([0,1],|\cdot|)$.  But Example 1 on \cite{3}, proves that $S\notin D(a,b)$ for $0\leq a,b\leq1$, proving in this way that in fact $S\notin D(a,b)$ when the mapping is defined in the cone metric space $E$.  On the other hand, notice that if $S\in D_T(a,b)$ then next inequality follows
\begin{equation*}
d(\alpha \sqrt{x},\alpha\sqrt{y})\leq ad(\alpha x,\alpha y)+b[d(\alpha x,\alpha\sqrt{x})+d(\alpha y,\alpha\sqrt{y})]
\end{equation*}
or, that $S\in D(a,b)$ which is false.
\end{ejem}

\begin{ejem}
Now, we will consider a function $S\notin D(a,b)$ and an injective, continuous function $T$ such that $S\in D_T(a,b)$.

In this case $E=(C[0,1],\R)$, the cone $P=\{\varphi\in E\;:\;\varphi\geq0\}$,  $M:=[\frac{1}{2},1]$ and $d(x,y)=|x-y|e^t$ where $e^t\in E$, consider the function $S$ define as $Sx=\sqrt{x}$ which does not belong to $D(a,b)$ for $0\leq a,b\leq1$. Let the continuous and injective function $Tx=\ln x$. We will prove that $S\in D_T(a,b)$. For attain a such goal we need to proof that the following inequality holds.
\begin{equation*}
|\ln\sqrt{x}-\ln\sqrt{y}|\leq a|\ln x-\ln y|+b[|\ln x-\ln\sqrt{x}|+|\ln y-\ln\sqrt{y}|].
\end{equation*}
Or equivalently that
\begin{equation}\label{ine:exam}
\left|\ln\sqrt\frac{x}{y}\right|\leq a\left|\ln\frac{x}{y}\right|+b\left[\left|\ln\frac{\sqrt{x}}{x}\right|+\left|\ln\frac{y}{\sqrt{y}}\right|\right].
\end{equation}
Notice that inequality above holds for $x=y$. So, we can assume without loss of generality that $x>y$. Then \eqref{ine:exam} ca be rewrite as
\begin{equation*}
\ln\sqrt\frac{x}{y}\leq a\ln\frac{x}{y}+b\ln\frac{\sqrt x}{x}-b\ln\frac{y}{\sqrt y}.
\end{equation*}    
Since $0\leq a,b<1$ and $x,y\in M$, then we have that last inequality is satisfies if and only if
\begin{equation*}
\sqrt\frac{x}{y}\leq e^{a+b}\frac{x}{y}\frac{\sqrt{x}}{x}\frac{y}{\sqrt{y}}
\end{equation*}
which is holds if and only if $a+b\geq0$. Therefore, $S\in D_T(a,b)$.
\end{ejem}

\begin{ejem}
In this example we will show a function $S\in D(a,b)$ such that for a continuous and injective function $T$, we get $S\notin D_T(a,b)$.

Let $E=(C[0,1],\R)$, the cone $P=\{\varphi\in E\;:\;\varphi\geq0\}$,  $M:=[0,1]$ and $d(x,y)=|x-y|e^t$ where $e^t\in E$. In this case let the function $Sx=x$. We have that $S\in D(a,b)$ for $a<1$. Take $Tx=\sqrt{x}$. Then if $S\in D_T(a,b)$, the following inequality is satisfied
\begin{equation*}
d(\sqrt x,\sqrt y)\leq ad(\sqrt x,\sqrt y),
\end{equation*}
rewriting inequality above we get
\begin{equation*}
|\sqrt x-\sqrt y|\leq a|\sqrt x-\sqrt y|,
\end{equation*}
which is false because $a<1$.
\end{ejem}

The following result describe the asymptotic behavior of the elements on $D_T(a,b)$, also give conditions for existence and uniqueness of fixed points for applications in that class.

\begin{teo}\label{teo3.3}
Let $(M,d)$ be a complete cone metric space, $P$ be a cone, and let $T,S: M\longrightarrow M$ be mappings such that $T$ is continuous and injective. If $S\in D_T(a,b)$, $a,b\geq0$ and moreover:
\begin{enumerate}
\item[(i)] If $a+2b<1$, then
\begin{equation*}
\tlim_{m,n\to\infty}d(TS^mx_0,TS^nx_0)=0,\; \mbox{for every $x_0\in X$.}
\end{equation*}
\item[(ii)] If $a+2b<1$, and $T$ is subsequentially convergent. Then $S$ has a fixed point.

\item[(iii)] If $a<1$, then the fixed point of $S$ is unique.
\item[(iv)] If $a+2b<1$ and furthermore $T$ is sequentially convergent, then the sequence of iterates $(S^nx_0)$, for every $x_0\in X$, converges to the unique fixed point of $S$.
\item[(v)] If $b<1$ and $(x_n)_n$ is a sequence of points on $X$ converging to a fixed point of $S$, then
\begin{equation}\label{eq:thPrinc1}
d(Tx_n,TSx_n)\to0,\quad\mbox{as $n\to\infty$.}
\end{equation}
\item[(vi)] If $b<1$, $(x_n)_n$ converges to a point $p$ and  \eqref{eq:thPrinc1} holds, then $p$ is a fixed point of $S$.
\item[(vii)]If $b<1$ and $S$ has a fixed point $p$, then $S$ is continuous at $p$.
\end{enumerate}
\end{teo}
\begin{proof}
\begin{itemize}
\item[(i)] Let $x_0$ be an arbitrary point, we define the iterative sequence $(x_n)_n$ by setting $x_{n+1}:=Sx_n$ (equivalently, $x_n:=S^nx_0$), $n=1,\dots$. Using the fact that $S\in D_T(a,b)$, we have
    \begin{eqnarray*}
      d(Tx_n,Tx_{n+1}) &=& d(TSx_{n-1},TSx_n) \\
       &\leq& ad(Tx_{n-1},Tx_n)+b[d(Tx_{n-1},TSx_{n-1})+d(Tx_n,TSx_n)]
    \end{eqnarray*}
    so that,
    \begin{equation*}
    d(Tx_n,Tx_{n+1})\leq\frac{a+b}{1-b}d(Tx_{n-1},Tx_n).
    \end{equation*}
    Repeating the argument we obtain
    \begin{equation*}
    d(Tx_n,Tx_{n+1})\leq\lambda^nd(Tx_0,Tx_1),
    \end{equation*}
    where $\lambda:=\frac{a+b}{1-b}<1$. From inequality above, we have for every $m,n\in\N$ such that $m>n$
    \begin{eqnarray*}
      d(Tx_m,Tx_n) &\leq&d(Tx_m,Tx_{m-1})+d(Tx_{m-1},Tx_{m-2})\\
      &&+\dots+d(Tx_{n+1},Tx_n)  \\
       &\leq& [\lambda^{m-1}+\lambda^{m-2}+\dots+\lambda^n]d(Tx_0,Tx_1) \\
       &\leq&  [\lambda^n+\lambda^{n-1}+\dots+]d(Tx_0,Tx_1)=\frac{\lambda^n}{1-\lambda}d(Tx_0,Tx_1).
    \end{eqnarray*}
    Let $0\ll c$ be given. Choose a natural number $N_1$ such that $\frac{\lambda^n}{1-\lambda}d(Tx_0,Tx_1)\ll c$, for all $n\geq N_1$. Thus
    \begin{equation*}
    d(Tx_m,Tx_n)\ll c,\qquad m>n.
    \end{equation*}
    So, we have that $d(Tx_m,Tx_n)\ll \frac{c}{n}$ for all $n\geq1$. Therefore, $\frac{c}{n}-d(Tx_m,Tx_n)\in P$ for all $n\geq1$. Since $\frac{c}{n}\to0$ ($n\to\infty$) and $P$ is closed, then $-d(Tx_m,Tx_n)\in P$. Thus
    \begin{equation*}
    \lim_{n,m\to\infty}d(TS^mx_0,TS^nx_0)=0.
    \end{equation*}

\item[(ii)] From part (i) we know that $(Tx_n)_n$ is a Cauchy sequence in $(X,d)$. Since $(X,d)$ is a complete cone metric space, there exists $x\in X$ such that
    \begin{equation}\label{eq:ThmPric2}
    \lim_{n\to\infty}Tx_n=x.
    \end{equation}
    On the other hand, the hypothesis that $T$ is a subsequentially convergent mapping, imply that  $(x_n)_n$ has a convergent subsequence. So, there exists $y\in X$ and $(x_{n(k)})_k$ such that $\lim_{k\to\infty}x_{n(k)}=y$. From the continuity of $T$ we have that $\lim_{k\to\infty}Tx_{n(k)}=Ty$ and by equality \eqref{eq:ThmPric2} we get $Ty=x$. Now, choose a natural number $N_2$ such that $d(Ty,Tx_{n(k)-1})\ll\frac{1-b}{4a}c$, $\;\lambda^{n(k)-1}d(TSx_0,Tx_0)\ll\frac{1-b}{4b}c$, $\;\lambda^{n(k)}d(TSx_0,Tx_0)\ll\frac{1-b}{4}c$ and $d(Tx_{n(k)+1},Ty)\ll\frac{1-b}{4}c$, for all $n(k)>N_2$.
    \begin{eqnarray*}
      d(TSy,Ty) &\leq& d(TSy,TS^{n(k)}x_0)+d(TS^{n(k)}x_0,TS^{n(k)+1}x_0)\\
      &&+d(TS^{n(k)+1}x_0,Ty) \\
       &\leq& ad(Ty,TS^{n(k)-1}x_0)+b[d(Ty,TSy)\\
       &&+d(TS^{n(k)-1}x_0,TS^{n(k)}x_0)]  \\
       &&+d(TS^{n(k)}x_0,TS^{n(k)+1}x_0)+d(TS^{n(k)+1}x_0,Ty)  \\
       &\leq&ad(Ty,TS^{n(k)-1}x_0)+b[d(Ty,TSy)\\
       &&+\lambda^{n(k)-1}d(TSx_0,Tx_0)]  \\
       &&+\lambda^{n(k)}d(TSx_0,Tx_0)+d(TS^{n(k)+1}x_0,Ty),
   \end{eqnarray*}
   hence, we get
   \begin{eqnarray*}
     d(TSy,Ty) &\leq& \frac{a}{1-b}d(Ty,Tx_{n(k)-1})+\frac{b}{1-b}\lambda^{n(k)-1}d(TSx_0,Tx_0) \\
      && +\frac{1}{1-b}\lambda^{n(k)}d(TSx_0,Tx_0)+\frac{1}{1-b}d(Tx_{n(k)+1},Ty)\\
      &\ll&\frac{c}{4}+\frac{c}{4}+\frac{c}{4}+\frac{c}{4}=c,
   \end{eqnarray*}
   for all $n(k)>N_2$. Thus $d(TSy,Ty)\ll\frac{c}{m}$ for all $m\geq1$. So, $\frac{c}{m}-d(TSy,Ty)\in P$ for all $m\geq1$. Since $\frac{c}{m}\to\infty$, as $m\to\infty$, and $P$ is closed, then $-d(TSy,Ty)\in P$. Therefore, $d(TSy,Ty)=0$. The injectivity of $T$ imply that $d(Sy,y)=0$, So $S$ has a fixed point.
\item[(iii)] Let $y_1$ and $y_2$ be two fixed points of $S$. Using the fact that $S\in D_T(a,b)$, the following inequality is satisfied
    \begin{equation*}
    d(Ty_1,Ty_2)=d(TSy_1,TSy_2)\leq ad(Ty_1,Ty_2)
    \end{equation*}
    which only holds of $Ty_1=Ty_2$, or since $T$ is injective, $y_1=y_2$.
\item[(iv)] If $T$ is sequentially convergent, then replacing the sequence $(n(k))_k$ by $(n)_n$ we conclude that $\lim_{n\to\infty}x_n=y$. Showing in this way that $(x_n)_n$ converges to the fixed point of $S$.
\item[(v)] Let $x_n\to p$, where $Sp=p$. Then
\begin{equation}\label{eq:ThmPric3}
d(TSx_n,TSp)\leq ad(Tx_n,Tp)+bd(Tx_n,TSx_n).
\end{equation}
Applying the triangle inequality to the left side of inequality above, we obtain
\begin{equation*}
0\leq(1-b)d(TSx_n,Tx_n)\leq(1+a)d(Tx_n,Tp).
\end{equation*}
Since $T$ is continuous, then for $0\ll c$ given we can choose $N_1\in\N$ such that
\begin{equation*}
d(Tx_n,Tp)\ll \frac{c}{1-b},
\end{equation*}
thus, $d(TSx_n,Tx_n)\ll c$. Therefore, $d(TSx_n,Tx_n)\ll\frac{c}{n}$ for all $n\geq1$. Then we can conclude that $\frac{c}{n}-d(TSx_n,Tx_n)\in P$ for all $n\geq1$. Since $\frac{c}{n}\to 0$, as $n\to\infty$, and $P$ is closed, we have that $-d(TSx_n,Tx_n)\in P$ which finally imply $\displaystyle\lim_{n\to\infty}d(TSx_n,Tx_n)=0$.
\item[(vi)] By the triangle inequality we have
\begin{equation*}
d(Tp,TSp)\leq d(Tp,Tx_n)+d(Tx_n,TSx_n)+d(TSx_n,TSp)
\end{equation*}
in view of that $S\in D_T(a,b)$, we get
\begin{eqnarray*}
  d(Tp,TSp) &\leq& d(Tp,Tx_n)+d(Tx_n,TSx_n)+ad(Tx_n,Tp) \\
   && +b[d(Tx_n,TSx_n)+d(Tp,TSp)],
\end{eqnarray*}
therefore,
\begin{equation*}
0\leq(1-b)d(Tp,TSp)\leq(1+a)d(Tp,Tx_n)+(1-b)d(Tx_n,TSx_n).
\end{equation*}
From \eqref{eq:thPrinc1}, reasoning as in part (v), and because $T$ is injective, we have that $Sp=p$.
\item[(vii)] This in an immediate consequence of part (v), equality \eqref{eq:ThmPric3}, the continuity and injectivity of $T$.
\end{itemize}
\end{proof}

\section{Fixed point theorems  that are consequence of Theorem \ref{teo3.3} }

This section is devoted to exhibit how Theorem \ref{teo3.3} extends  the  fixed point theorems given in \cite{BMOP,KK,9,KhSa08,10,Mo09,MoRo09,RH}. We would like to put out that the mentioned Theorem \ref{teo3.3} also generalize the asymptotic results of the cited references.

\subsection*{Fixed point theorems on complete metric spaces}

If we take $E=\R_+:=\{x\in\R\;:\;x\geq0\}$ and $Tx=x$ in Theorem \ref{teo3.3}, then we obtain (for the corresponding value of $a$ and $b$) the following:
\begin{coro}[E.g. \cite{KK}, Theorem 2.1]
Let $(M,d)$ be a complete metric space and $S: M\longrightarrow M$ is a contraction mapping. Then $S$ has a unique fixed point.
\end{coro}

\begin{coro}[E.g. \cite{9}, Theorem 1]
Let $(M,d)$ be a complete metric space and $S: M\longrightarrow M$  a mapping satisfying,
\begin{equation*}
d(Sx,Sy)\leq b[d(x,Sx)+d(y,Sy)],
\end{equation*}
for all $x,y\in M$ and $b\in[0,\frac{1}{2})$. Then $S$ has a unique fixed point.
\end{coro}

\begin{coro}[\cite{10}, Theorem 2]
Let $(M,d)$ be a complete metric space, and $S:M\longrightarrow M$ an arbitrary mapping. If $S$
satisfies the following condition, for $a,b\geq0$ with $a+2b<1$ and any $x,y\in M$,
\begin{equation*}
d(Sx,Sy)\leq ad(x,y)+b\left[d(x,Sx)+d(y,Sy)\right].
\end{equation*}
Then $S$ has a unique fixed point.
\end{coro}

Now, if we only take $E=\R_{+}$ in the Theorem \ref{teo3.3} we obtain:

\begin{coro}[\cite{BMOP}, Theorem 2.6]
Let $(M,d)$ be a complete metric space and $T: M\longrightarrow M$ be an injective, continuous and subsequentially convergent mapping. Then every  continuous function $S: M\longrightarrow M$  satisfying
\begin{equation*}
d(TSx,TSy)\leq ad(Tx,Ty),
\end{equation*}
for all $x,y\in M$ and $0\leq a<1$, has a unique fixed point. Moreover, if $T$ is sequentially convergent, then for each $x_0\in M,$ the sequence $(S^nx_0)$ converge to the fixed point of $S$.
\end{coro}

\begin{coro}[\cite{Mo09}, Theorem 2.1]
Let $(M,d)$ be a complete metric space and $T,S: M\longrightarrow M$ be mappings such that $T$ is continuous, injective and
 subsequentially convergent. If $b\in[0,\frac{1}{2})$ and
 \begin{equation*}
d(TSx,TSy)\leq b\left[d(Tx,TSx)+d(Ty,TSy)\right],
\end{equation*}
  for all $x,y\in M$, then $S$ has a unique fixed point. Also if $T$ is sequentially convergent, then for each $x_0\in M,$ the sequence $(S^nx_0)$ converge to the fixed point of $S$.
\end{coro}

Next result is the version of $D_T(a,b)$ on metric spaces.
\begin{coro}\label{cor: DT(a,b) metric}
Let $(M,d)$ be a complete metric space and $T,S: M\longrightarrow M$ be mappings such that $T$ is continuous, injective and
 subsequentially convergent. If $a,b\geq0$, with $a+2b<1$ and
 \begin{equation*}
d(TSx,TSy)\leq ad(Tx,Ty)+b\left[d(Tx,TSx)+d(Ty,TSy)\right],
\end{equation*}
  for all $x,y\in M$, then $S$ has a unique fixed point. Also if $T$ is sequentially convergent, then for each $x_0\in M,$ the sequence $(S^nx_0)$ converge to the fixed point of $S$.
\end{coro}

\subsection*{Fixed point theorems on complete cone metric spaces}

As in previous part, for the corresponding different values of $a$ and $b$, if we assume $Tx=x$ in Theorem \ref{teo3.3} we get:

\begin{coro}[\cite{RH}, Theorem 2.3]
Let $(M,d)$ be a complete cone metric space and the mapping $S: M\longrightarrow M$ is a contraction function. Then $S$ has a unique fixed point in $M$ and for any $x_0\in M$, we have that $(S^nx_0)$ converges to the fixed point.
\end{coro}

\begin{coro}[\cite{RH}, Theorem 2.6]
Let $(M,d)$ be a complete cone metric space and the mapping $S: M\longrightarrow M$ satisfying,
\begin{equation*}
d(Sx,Sy)\leq b[d(x,Sx)+d(y,Sy)],
\end{equation*}
for all $x,y\in M$ and $b\in[0,\frac{1}{2})$. Then $S$ has a unique fixed point and for any $x_0\in M$, we have that $(S^nx_0)$ converges to the fixed point.
\end{coro}

\begin{coro}[C.f. \cite{KhSa08}, Theorem 2.7]
Let $(M,d)$ be a complete cone metric space and $S:M\longrightarrow M$ an arbitrary mapping. If $S$
satisfies the following condition, for $a,b\geq0$ with $a+2b<1$ and any $x,y\in M$,
\begin{equation*}
d(Sx,Sy)\leq ad(x,y)+b\left[d(x,Sx)+d(y,Sy)\right].
\end{equation*}
Then $S$ has a unique fixed point.
\end{coro}
The final case is to consider $T\neq \mbox{id}$, where by $\mbox{id}$ we denote the identity function $\mbox{id}x=x$ for all $x\in M$.

\begin{coro}[\cite{MoRo09}, Theorem 3.1]
Let $(M,d)$ be a complete cone metric space, $T: M\longrightarrow M$ be an one to one and continuous function, moreover let $S: M\longrightarrow M$ a  continuous function satisfying
\begin{equation*}
d(TSx,TSy)\leq ad(Tx,Ty),
\end{equation*}
for all $x,y\in M$ and $0\leq a<1$. Then
\begin{enumerate}
\item[(i)] For every $x_0\in M$,
\begin{equation*}
\tlim_{n\rightarrow \infty} d(TS^nx_0,TS^{n+1}x_0)=0;
\end{equation*}
\item[(ii)] There is $y_0\in M$ such that
\begin{equation*}
\tlim_{n\rightarrow \infty} TS^nx_0=y_0;
\end{equation*}
\item[(iii)] If $T$ is subsequentially convergent, then $(S^nx_0)$ has a convergent subsequence;
\item[(iv)] There is a unique $z_0\in M$ such that
\begin{equation*}
Sz_0=z_0;
\end{equation*}
\item[(v)] If $T$ is sequentially convergent, then for each $x_0\in M$ the iterate sequence $(S^nx_0)$
      converges to $z_0.$
\end{enumerate}
\end{coro}
The following result generalize the $T$-Kannan class given in  \cite{Mo09}.
\begin{coro}\label{cor:kannan cone}
Let $(M,d)$ be a complete cone metric space and $T,S: M\longrightarrow M$ be mappings such that $T$ is continuous, injective and
 subsequentially convergent. If $b\in[0,\frac{1}{2})$, and
 \begin{equation*}
d(TSx,TSy)\leq b\left[d(Tx,TSx)+d(Ty,TSy)\right],
\end{equation*}
  for all $x,y\in M$, then $S$ has a unique fixed point. Also if $T$ is sequentially convergent, then for each $x_0\in M,$ the sequence $(S^nx_0)$ converge to the fixed point of $S$.
\end{coro}

Finally, we would like summarize these previous results in the following table: Notice that each box in the table below indicate the reference of the theorem generalized by Theorem \ref{teo3.3} for each case of the function $T$ in the left upper corner.

\begin{table}[H]
\centering
\begin{tabular}{|c|c|c|c|}
  \hline
 $\begin{array}{c}
   T=\mbox{id}\\
  T\neq\mbox{id}\end{array}$& $a,b\geq0$ & $\begin{array}{c}
                      a=0 \\
                      b\in[0,\frac{1}{2})
                    \end{array}
  $ & $\begin{array}{c}
                      b=0 \\
                      a<1
                    \end{array}
  $  \\ \hline
  $\begin{array}{c}
     (M,d) \mbox{ be a complete}\\
     \mbox{cone metric space}
   \end{array}
  $  &  $\begin{array}{c}
   \mbox{C.f., \cite[Thm 2.7]{KhSa08}}\\
   D_T(a,b)
 \end{array}$ & $\begin{array}{c}
   \mbox{\cite[Thm 2.6]{RH}}\\
   D_T(0,b)
 \end{array}$ & $\begin{array}{c}
  \mbox{\cite[Thm 2.3]{RH}}\\
  \mbox{\cite[Thm 3.1]{MoRo09}}
 \end{array}$ \\ \hline
  $E=\R$ & $\begin{array}{c}
   \mbox{\cite[Thm 2]{10}}\\
   D_T(a,b)
 \end{array}$ & $\begin{array}{c}
   \mbox{\cite[Thm 1]{9}}\\
  \mbox{\cite[Thm 2.1]{Mo09}}
 \end{array}$ & $\begin{array}{c}
   \mbox{\cite[Thm 2.1]{KK} }\\
   \mbox{\cite[Thm 2.6]{BMOP}}
 \end{array}$ \\ \hline
\end{tabular}
\vspace{0.5cm}
 \caption{Results generalized by Theorem \ref{teo3.3}.}
\end{table}

\end{document}